# Limits of Mellin coefficients and Berezin transform

Benoit Barusseau

I.M.B, Université Bordeaux 1
351 cours de la Libération, 33400 Talence
*Email:* benoit.barusseau@math.u-bordeaux1.fr

**Abstract**

We consider a bounded function $f$ on [0,1). In [4] B. Korenblum and K.E.Zhu give a case where we have equality between the limit near the boundary of the unit disc of the Berezin transform and the limit of the normalized Mellin coefficient when one of them is 0. We previously describe the case where one limit point has modulus $\|f\|_\infty$. We also use the mean values of $f$ near 1.

The aim of this article is to show that between these two extreme cases, we can have distinct limits. In the same time, we also study the sets of limit points of these quantities.

## 1 Preliminaries

In the following, $\mathbb{D}$ denote the open unit disc and dA the normalized Lebesgue measure of $\mathbb{D}$.

The bergman space $L_a^2$ of the disc is the Hilbert space of analytics squared integrable functions on $\mathbb{D}$. A function $f \in L^\infty(\mathbb{D}, \mathrm{d}A)$ is said to be radial if $f(z) = f(|z|)$ almost everywhere on $\mathbb{D}$. The Toeplitz operator $T_f: L_a^2 \longrightarrow L_a^2$ is defined to be $T_f(h) = P(fh)$ where is the orthogonal projection from $L^2(\mathbb{D}, \mathrm{d}A)$ onto $L_a^2$.

Let $f$ be a bounded radial function, we denote $C_n(f) = (n+1)\int_0^1 f(r)dr$ the normalized Mellin coefficient of indice $n \in \mathbb{N}$. The Berezin transform of $f$ is defined to be

$$\tilde{f}(z) = \int_\mathbb{D} f(w)\frac{(1-|z|^2)^2}{|1-\bar{z}w|^4}\mathrm{d}A(w).$$

In [1], we show that

$$0 \leqslant \limsup_{z \to \partial\mathbb{D}} |\tilde{f}(z)| \leqslant \limsup_{n \to \infty} |C_n(f)| = \|T_f\|_e \leqslant \limsup_{\varepsilon \to 1} \left|\frac{1}{1-\varepsilon}\int_\varepsilon^1 f(r)dr\right| \leqslant \|f\|_\infty. \tag{1}$$

where $\|T_f\|$ denote the essential norm of the bounded Toeplitz operator on the Bergman space.

Now considering the three intermediate quantities, B.Korenblum and K.E.Zhu show in [4] that, if one of these quantities is zero so is for the other. We have shown in [1] that if one of the previous quantity is equal to $\|f\|_\infty$ then so is for the other. Moreover, these two works answer the questions: when the essential norm of the Toeplitz operator $T_f$ is zero and when it is $\|f\|_\infty$.

Here, we answer the following natural question: is there exist a bounded function which gives strict inequalities in equation (1)?

We prove that such a function exists and explain why it works. Our method underlines a link between the limits points of the normalized Mellin coefficients and the limits of the Berezin transform near the boundary of the unit disc.

We will need some technical tools which are contained in the following lemmas.

**Lemma 1.** *(Ž. Čučković [2])Let $F$ be a bounded m-quasihomogeneous function with $f$ as radial part. For any complex $z = Re^{i\theta}$, we have*

$$\tilde{F}(z) = 2(1-R^2)^2 R^{|m|}e^{im\theta}\sum_{n=0}^\infty \frac{n(n+|m|)}{2n+|m|+1}C_{2n+|m|}(f)R^{2(n-1)}.$$

**Lemma 2.** *(B. Barusseau [1]) Let $\rho > 0$, $K$ a compact subset of $\overline{\rho\mathbb{D}}$ and $L \in \rho\mathbb{T}\backslash K$. If $(a_R)_{R\in[0,1)}$, $(M_R)_{R\in[0,1)}$ and $(N_R)_{R\in[0,1)}$ are complex functions of $R \in [0,1)$ such that*

$$a_R \in [0,1], M_R \in K, N_R \in \overline{\rho\mathbb{D}}$$

*then*

$$\liminf_{R\to 1} a_R > 0 \Longrightarrow \liminf_{R\to 1} |a_R M_R + (1-a_R)N_R - L| > 0.$$

## 2 The result

**Theorem 3.** *There exists $f$ a bounded radial function such that the sets*

1. $A_M(f) = \{L \in \mathbb{C}, \liminf_{z\to\partial\mathbb{D}} |L - C_n(f)|\}$ *of the limits points of the normalized Mellin coefficients,*

2. $A_{\text{me}}(f) = \left\{L \in \mathbb{C}, \liminf_{\varepsilon\to 1} \left|L - \frac{1}{1-\varepsilon}\int_\varepsilon^1 f(t)dt\right|\right\}$ *of the limits points of the mean values of $f$ near 1,*

3. $A_B(f) = \left\{L \in \mathbb{C}, \liminf_{z\to\partial\mathbb{D}} \left|L - \tilde{f}(z)\right|\right\}$ *of the limits points of the Berezin transform near the boundary of the unit disc*

*are disjoint. Moreover, these three sets have disctinct supremum, this is*

$$\limsup_{z\to\partial\mathbb{D}} |\tilde{f}(z)| < \limsup_{n\to\infty} |C_n(f)| < \limsup_{\varepsilon\to 1} \left|\frac{1}{1-\varepsilon}\int_\varepsilon^1 f(r)dr\right|.$$

Using proposition 3 in [1], we know that if one of these quantities has a limit (not only a superior limit) then so is for the others and these limits are equal. Thus we need to consider a function which has no limit in the neighbourhood of 1. A simple way is to look at functions which quickly oscillate in the neighbourhood of 1. S.Grudsky and N. Vasilevski define in [3] the following function $f$: for any $z \in \mathbb{D}$

$$f(z) = \alpha \ln(|z|^{-1})^i \tag{2}$$

where $\alpha = \left(\int_0^1 \ln(r^{-1})^i dr\right)^{-1}$. From now, $f$ is fixed by equation (2). We remark that $f$ is bounded but has no limit in 1.

**Proposition 4.** *[3] We have*
$$\{C_n(f), n \in \mathbb{N}\} = \{\exp(-i\ln(n+1)), n \in \mathbb{N}\}$$
*In particular $A_M(f) = \partial\mathbb{D}$ and $\limsup_{n\to\infty} |C_n(f)| = \|T_f\|_e = 1$.*

Now, we consider the set $A_{\text{me}}(f)$.

**Proposition 5.** *We have*
$$A_{\text{me}}(f) = \frac{|\alpha|}{\sqrt{2}}\partial\mathbb{D}.$$
*In particular, $\limsup_{\varepsilon\to 1}\left|\frac{1}{1-\varepsilon}\int_\varepsilon^1 f(r)dr\right| = \frac{|\alpha|}{\sqrt{2}} > 1$.*

**Proof.**

Let $\varepsilon > 0$, a change of variable gives
$$\frac{1}{1-\varepsilon}\int_\varepsilon^1 f(r)dr = \frac{\alpha}{1-\varepsilon}\int_\varepsilon^1 \ln(1/r)^i dr$$
$$= \frac{\alpha}{1-\varepsilon}\int_{-\infty}^{\ln\ln(1/\varepsilon)} e^{iv}e^v e^{-e^v} dv$$

We split the previous integrale in the following way :
$$\frac{1}{1-\varepsilon}\int_\varepsilon^1 f(r)dr = \frac{\alpha}{1-\varepsilon}\int_{-\infty}^{\ln\ln(1/\varepsilon)} e^{iv}e^v dv + \beta(\varepsilon) \tag{3}$$

and $\beta(\varepsilon) = \frac{\alpha}{1-\varepsilon} \int_{-\infty}^{\ln\ln(1/\varepsilon)} e^{iv} e^v (e^{-e^v} - 1) dv$. Using the mean value theorem, we obtain that for any $u \in \mathbb{R}$,

$$\left| e^{-e^u} - 1 \right| \leqslant \sup_{v<0} |e^v| |-e^u - 0| \leqslant e^u.$$

thus

$$|\beta(\varepsilon)| \leqslant \frac{|\alpha|}{1-\varepsilon} \int_{-\infty}^{\ln\ln(1/\varepsilon)} \left| e^{(i+2)v} \right| dv$$

$$\leqslant \frac{|\alpha|}{1-\varepsilon} \int_{-\infty}^{\ln\ln(1/\varepsilon)} e^{2v} dv = \frac{\alpha}{1-\varepsilon} \frac{1}{2} \ln\left(\frac{1}{\varepsilon}\right)^2.$$

Now, since

$$\lim_{\varepsilon \to 1} \frac{\ln(1/\varepsilon)}{1-\varepsilon} = \lim_{\varepsilon \to 1} \frac{1}{\varepsilon} = 1,$$

we have

$$\lim_{\varepsilon \to 1} |\beta(\varepsilon)| \leqslant \frac{|\alpha|}{2} \lim_{\varepsilon \to 1} \ln(1/\varepsilon) = 0.$$

Equation (3) gives

$$\left| \frac{1}{1-\varepsilon} \int_{\varepsilon}^{1} f(r) dr \right| = \left| \frac{\alpha}{1-\varepsilon} \int_{-\infty}^{\ln\ln(1/\varepsilon)} e^{iv} e^v dv + o_{\varepsilon \to 1}(1) \right|$$

$$= \left| \frac{\alpha}{1-\varepsilon} \frac{1}{1+i} \ln(1/\varepsilon)^{i+1} + o_{\varepsilon \to 1}(1) \right|$$

$$= \left| \frac{\alpha}{1+i} \ln(1/\varepsilon)^i \frac{\ln(1/\varepsilon)}{1-\varepsilon} + o_{\varepsilon \to 1}(1) \right|$$

Since $\lim_{\varepsilon \to 1} \left| \ln(1/\varepsilon)^i \frac{\ln(1/\varepsilon)}{1-\varepsilon} \right| = 1$, we have

$$\lim_{\varepsilon \to 1} \left| \frac{1}{1-\varepsilon} \int_{\varepsilon}^{1} f(r) dr \right| = \left| \frac{\alpha}{1+i} \right| = \frac{|\alpha|}{\sqrt{2}}.$$

It is clear that $A_{\mathrm{me}}(f) \subset \frac{|\alpha|}{\sqrt{2}} \partial \mathbb{D}$. The converse is clearly true. Finally, a simple estimate of $|\alpha|$ gives that $\frac{|\alpha|}{\sqrt{2}} > 1$. $\square$

To describe $A_M(f)$, we will use lemma 2 and summation by part. We recall that considering two sequences $(a_n)_{n \in \mathbb{N}}$ and $(b_n)_{n \in \mathbb{N}}$, for any $N \in \mathbb{N}$, we have

$$\sum_{n=0}^{N} a_n b_n = a_N \sum_{n=0}^{N} b_n - \sum_{n=0}^{N-1} \left( \sum_{k=0}^{n} b_k \right) (a_{n+1} - a_n).$$

**Theorem 6.** *Let $g$ be a bounded radial function and $L$ be such that $|L| = \sup A_M(g)$. If there exists $\varepsilon > 0$ such that*

$$\liminf_{N \to \infty} \frac{\mathrm{Card}\{n \in \{0, ..., N\}, |C_{2n+1}(g) - L| > \varepsilon\}}{N} > 0 \tag{4}$$

*then $L \notin A_B(g)$.*

*In particular, if it is true for any $L \in \sup A_M(g) \partial \mathbb{D}$ then $\|T_f\|_e > \limsup_{z \to \partial \mathbb{D}} |\tilde{f}(z)|$.*

**Proof.** To show this proposition, it suffices to show that

$$\liminf_{z \to \partial \mathbb{D}} |\tilde{f}(z) - L| > 0. \tag{5}$$

In order to simplify our calculations, we consider the 1-quasihomogeneous function $F$ with radial part $f$. This means that

$$F(re^{i\theta}) = e^{i\theta} f(r).$$

We show in [1] that equation (5) is equivalent to

$$\liminf_{z \to \partial \mathbb{D}} |\tilde{F}(z) - L| > 0.$$

Since $\tilde{F}$ is also 1-quasihomogeneous, the previous equality is obviously equivalent to the following one

$$\liminf_{R \to 1} |\tilde{F}(R)/R - L| > 0$$

We will apply lemma 2 with $\rho = |L| = \sup A_M(g)$ and

$$K = \overline{\mathcal{EC}(\mathbb{D} \setminus D(L, \varepsilon))}$$

the closure of the convex hull of $\mathbb{D} \setminus D(L, \varepsilon)$.

We denote

$$P = \{n \in \mathbb{N}, |C_{2n+1}(g) - L| > \varepsilon\}$$

and for all $0 < R < 1$,

$$M_R = \frac{\sum_{n \in P} n C_{2n+1}(g) R^{2n-2}}{\sum_{n \in P} n R^{2n-2}},$$

$$N_R = \frac{\sum_{n \in \mathbb{N} \setminus P} n C_{2n+1}(g) R^{2n-2}}{\sum_{n \in \mathbb{N} \setminus P} n R^{2n-2}}$$

and

$$a_R = (1 - R^2)^2 \sum_{n \in P} n R^{2n-2}.$$

Thus, for all $0 < R < 1$, we have

$$a_R M_R + (1 - a_R) N_R = (1 - R^2)^2 \sum_{n \in \mathbb{N}} n C_n(g) R^{2n-2}$$

and

$$\tilde{F}(R)/R - L = a_R M_R + (1 - a_R) N_R - L.$$

Since for all $0 < R < 1$, we have

$$(1 - R^2)^2 \sum_{n \in \mathbb{N}} n R^{2n-2} = 1,$$

we see that $a_R \in [0, 1]$. Moreover, it is clear that $M_R \in K$ and $N_R \in \rho \overline{\mathbb{D}}$. Now in order to apply lemma 2, it remains to show that $\liminf_{R \to 1} a_R > 0$.

For any $n \in \mathbb{N}$, we denote $p_n = \text{Card}\{n \in \{0, ..., N\}, |C_{2n+1}(g) - L| > \varepsilon\}$.

⋆ Let $N \in \mathbb{N}^*$. Using a summation by part with $a_n = \begin{cases} n & \text{si } n \in P_N \\ 0 & \text{sinon} \end{cases}$ and $b_n = R^{2n-2}$, we have

$$\sum_{n \in P_N} n R^{2n-1} = R^{2N-2} \sum_{k \in P_N} k + \sum_{n=1}^{N-1} \left( \sum_{k \in P_n} k \right) (R^{2n-2} - R^{2n}). \tag{6}$$

We easily see that

$$\sum_{k \in P_n} k \geqslant \sum_{1 \leqslant k \leqslant p_n} k = \frac{p_n(p_n+1)}{2} \geqslant \frac{p_n^2}{2}.$$

Thus

$$\sum_{n=1}^{N-1} \left( \sum_{k \in P_n} k \right) (R^{2n-2} - R^{2n}) \geqslant \sum_{n=1}^{N-1} \frac{p_n^2}{2} (R^{2n-2} - R^{2n}).$$

By hypothesis there exists $M > 0$ such that $p_n > \delta n$ for $n > M$, so

$$\sum_{n=1}^{N-1} \frac{p_n^2}{2} (R^{2n-2} - R^{2n}) \geqslant \frac{\delta^2}{2} (1 - R^2) \sum_{n=M}^{N-1} n^2 R^{2n-2} \ .$$

Thus using equation (6), we deduce

$$\sum_{n \in P} n\, R^{2n-2} \geqslant \lim_{N \to \infty} \left( R^{2N-1} \sum_{k \in P_N} k \right) + \frac{\delta^2}{2}(1-R^2) \sum_{n=M}^{\infty} n^2 R^{2n-2}$$

$$\geqslant \frac{\delta^2}{2}(1-R^2) \sum_{n=M}^{\infty} n^2 R^{2n-2}.$$

Since

$$\sum_{n=0}^{\infty} n^2 R^{2n-2} = \frac{(R^2+1)}{(1-R^2)^3}.$$

we have

$$\lim_{R \to 1}\, (1-R^2) \sum_{n=M}^{\infty} n^2 R^{2n-2} = 2.$$

Thus

$$\liminf_{R \to 1} a_R \geqslant \delta^2 > 0.$$

Since $\partial \mathbb{D}$ is compact and $A_B(f)$ is closed, we have the particular case. $\square$

**Corollary 7.** *There exists $\gamma < 1$ such that*

$$A_B(f) \subset \gamma \mathbb{D}$$

*In particular, we have* $\limsup_{z \to \partial \mathbb{D}} |\tilde{f}(z)| < 1.$

**Proof.** Let $L \in \partial \mathbb{D}$, there exists $\theta \in\, ]-\pi, \pi]$ such that $L = e^{-i\theta}$. We will apply theorem 6 with $\varepsilon = \sqrt{2}$.

We recall that for any $N \in \mathbb{N}$,

$$n \in P_N \iff \left( 0 \leqslant n \leqslant N \text{ and } |e^{i\ln(2n+1)} - e^{i\theta}| > \sqrt{2} \right)$$

$$\iff \left( 0 \leqslant n \leqslant N \text{ and } \ln(2n+1) \in \bigcup_{k \geqslant -1} \left] \theta + \frac{\pi}{2} + 2k\pi, \theta + \frac{3\pi}{2} + 2k\pi \right[ \right)$$

For any sufficiently large $N \in \mathbb{N}$, there exists $t \in \mathbb{R}^+$ such that $N = \frac{1}{2}(e^{\theta + \frac{3\pi}{2} + 2t\pi} - 1)$. Thus

$$p_N \geqslant \sum_{k=0}^{[t]} \mathrm{Card}\left( \left] \frac{1}{2}(e^{\theta + \frac{\pi}{2} + 2k\pi} - 1), \frac{1}{2}(e^{\theta + \frac{3\pi}{2} + 2k\pi} - 1) \right[ \cap \mathbb{N} \right)$$

$$\geqslant \sum_{k=0}^{[t]} \left( \frac{1}{2} e^{\theta + \frac{3\pi}{2} + 2k\pi} - \frac{1}{2} e^{\theta + \frac{\pi}{2} + 2k\pi} - 1 \right)$$

Finally, we remark that

$$\sum_{k=0}^{[t]} \frac{1}{2} e^{\theta + \frac{3\pi}{2} + 2k\pi} - \frac{1}{2} e^{\theta + \frac{\pi}{2} + 2k\pi} - 1 \geqslant \sum_{k=0}^{[t]} \frac{1}{2} e^{\pi + \theta + 2k\pi} \geqslant K\left( e^{2([t]+1)\pi} - 1 \right)$$

where $K > 0$. Thus

$$p_N/N \geqslant K \frac{e^{2([t]+1)\pi} - 1}{\frac{1}{2} e^{\theta + \frac{3\pi}{2} + 2t\pi} - 1}$$

$$= K \frac{X - 1}{\frac{1}{2} e^{\theta + \frac{3\pi}{2}} X - 1}$$

where $X = e^{2([t]+1)\pi}$. It is now clear that $\liminf_N p_N/N > 0$ and using the previous theorem, the proof is finished. $\square$

Theorem 3 is now proved.

We have given a function $f$ such that $A_M(f)$, $A_B(f)$ and $A_{\mathrm{me}}(f)$ are disjoint sets and considering the real part of $f$, we obtain $A_B(\operatorname{Re} f) \subset A_M(\operatorname{Re} f) \subset A_{\mathrm{me}}(\operatorname{Re} f)$.

## 3 Generalization to extreme points

We can generalize this approach and give conditions on the sequence of the normalized Mellin coefficients to ensure that some particular points cannot be attained by both Berezin transform and normalized Mellin coefficients.

**Theorem 8.** *Let $g$ be a bounded radial function and $L$ be an extreme point of $A_M(g)$. If there exist $\varepsilon > 0$ such that*

$$\liminf_{N \to \infty} \frac{\operatorname{Card}\{n \in \{0, ..., N\}, |C_{2n+1}(g) - L| > \varepsilon\}}{N} > 0 \qquad (7)$$

*then $L \notin A_B$.*

**Proof.** Since $g$ is bounded so is $A_M(g)$. If $L$ is an extreme point of the closed set $A_M(g)$ then there exists $a \in \mathbb{C}$ and $r > 0$ such that $A_M \subset \overline{D(a, r)}$ and $A_M \cap \overline{D(a, r)} = L$.

Thus $\|g - a\|_\infty = |L - a| = r$. The previous corollary can be applied with $g - a$. And since

$$C_{2n+1}(g - a) - (L - a) = C_{2n+1}(g) - L$$

the condition is just condition (7) and we can conclude that

$$L - a \notin \left\{ L \in \mathbb{C}, \liminf_{z \to \partial \mathbb{D}} \left|(L - a) - \widetilde{g - a}(z)\right| \right\},$$

which is just $L \notin A_B$ since $\widetilde{g - a} = \tilde{g} - a$. $\square$

We can hope that if $A_M$ contains two elements then the condition (7) is true. This will be true if $A_M$ has a finite cardinal. In fact, since $A_M$ is a connected set (see [3]), if $A_M$ contains at least two elements then $A_M$ is not finite.

Using the same idea as in theorem 8, we give now an analogous of theorem 4 in [1] considering extreme points.

**Theorem 9.** *Let $g$ be a bounded radial function and $L$ be an extreme point of the essential range of $g$. The following conditions are equivalent*

*a) $L \in A_B(g)$;*

*b) $L \in A_M(g)$;*

*c) $L \in A_{\mathrm{me}}(g)$.*

**Proof.** Using the theorem 4 in [1] together with the same argument as in the proof of the previous theorem and the fact that $\frac{1}{1-\varepsilon} \int_\varepsilon^1 (g(r) - L) dr = \frac{1}{1-\varepsilon} \int_\varepsilon^1 g(r) dr - L$, we clearly obtain the theorem. $\square$

The following example is an application of the previous theorem.

**Example 10.** Let $g$ be defined by

$$\forall r \in [0, 1], g(r) = \begin{cases} 16r, r \in [0, 1/16[ \\ 1, r \in \bigcup_{n \geqslant 1} \left[1 - \frac{1}{2^{(2n)^2}}, 1 - \frac{1}{2^{(2n+1)^2}}\right[ \\ 0, r \in \bigcup_{n \geqslant 1} \left[1 - \frac{1}{2^{(2n+1)^2}}, 1 - \frac{1}{2^{(2n+2)^2}}\right[ \end{cases}.$$

For any $n \geqslant 1$, we have

$$\begin{aligned}
2^{(2n)^2} \int_{1-\frac{1}{2^{(2n)^2}}}^{1} g(r)dr &= 2^{(2n)^2} \sum_{p \geqslant n} \mu\left(\left[1-\frac{1}{2^{(2p)^2}}, 1-\frac{1}{2^{(2p+1)^2}}\right[\right) \\
&\geqslant 2^{(2n)^2} \mu\left(\left[1-\frac{1}{2^{(2n)^2}}, 1-\frac{1}{2^{(2n+1)^2}}\right[\right) \\
&\geqslant 2^{(2n)^2}\left(\frac{1}{2^{(2n)^2}} - \frac{1}{2^{(2n+1)^2}}\right) = 1 - \frac{1}{2^{4n+1}}.
\end{aligned}$$

Since $2^{(2n)^2} \int_{1-\frac{1}{2^{(2n)^2}}}^{1} g(r)dr \leqslant 1$, we have

$$\lim_{n \to \infty} 2^{(2n)^2} \int_{1-\frac{1}{2^{(2n)^2}}}^{1} g(r)dr = 1.$$

The same reasoning gives

$$2^{(2n+1)^2} \int_{1-\frac{1}{2^{(2n+1)^2}}}^{1} 1 - g(r)dr \geqslant 1 - \frac{1}{2^{6n+3}}.$$

Which implies that

$$0 \leqslant 2^{(2n+1)^2} \int_{1-\frac{1}{2^{(2n+1)^2}}}^{1} g(r)dr \leqslant \frac{1}{2^{6n+3}}$$

so

$$\lim_{n \to \infty} 2^{(2n+1)^2} \int_{1-\frac{1}{2^{(2n+1)^2}}}^{1} g(r)dr = 0.$$

Since $0 = \inf f$ and $1 = \sup f$ are extreme points of the essential range of $g$, we have $0, 1 \in A_M$ and $0, 1 \in A_B$.

We have given condition on the set of the Mellin coefficients to ensure that these limits are distinct or equal. It remains to give conditions on $f$ or on the means of $f$ near 1 in the case where they are distinct. It is still an open question.